\documentclass[a4paper]{amsart} 

\usepackage{amsfonts}
\usepackage{amsmath}  
\usepackage{amscd}
\usepackage{graphicx}
\usepackage{psfrag}

\usepackage{amssymb}
\usepackage{latexsym}
\usepackage{framed, color}

\newtheorem{thm}{\sc Theorem}[section]
\newtheorem{cor}[thm]{\sc Corollary}
\newtheorem{lem}[thm]{\sc Lemma}
\newtheorem{prop}[thm]{\sc Proposition}
\theoremstyle{definition}
\newtheorem{dfn}[thm]{\sc Definition}

\theoremstyle{remark}
\newtheorem{exam}[thm]{\sc Example}
\newtheorem{rmk}[thm]{\sc Remark}
\newtheorem{prob}[thm]{\sc Problem}
\makeatletter
 
 \@addtoreset{equation}{section}
\makeatother

\newcommand{\R}{\mathbf{R}}

\newcommand{\Z}{\mathbf{Z}}

\newcommand{\Int}{\mathop{\mathrm{Int}}\nolimits}

\newcommand{\rank}{\mathop{\mathrm{rank}}\nolimits}

\newcommand{\id}{\mathop{\mathrm{id}}\nolimits}

\renewcommand{\setminus}{\smallsetminus}
\def\spmapright#1{\smash{%
 \mathop{\hbox to 1.3cm{\rightarrowfill}}
  \limits^{#1}}}

\title[Special generic maps and Gromoll filtration]{Special
generic maps and Gromoll filtration}
\author{Osamu Saeki} 
\address{Institute of Mathematics for Industry,
Kyushu University,
Motooka 744, Nishi-ku, Fukuoka 819-0395, Japan}

\date{\today}

\begin{document}
\begin{abstract}
A smooth map of a closed $n$--dimensional manifold
into $\R^p$ with $1 \leq p \leq n$ is a special generic map
if it has only definite folds as its singularities. 
We show that for $1 \leq p < n$ and $n \geq 6$,
a homotopy $n$--sphere admits a special generic map into $\R^p$ 
with standard properties if and only if its Gromoll filtration
is equal to $p$.
\end{abstract}

\maketitle 

\section{Introduction}\label{section1}

Let $M$ be a smooth closed $n$--dimensional manifold.
In order to study the differential topology of $M$,
it is known to be quite beneficial to use
Morse functions on it.
For example, the homological information, or even
the information on the handlebody decomposition of $M$ can
be obtained by such a method. However, it is often quite
difficult to get the information on the attaching maps
of handles. For example, if there exists a Morse
function $f : M \to \R$ with exactly two critical
points, then $M$ is homeomorphic to the $n$--sphere $S^n$
\cite{Re1}: however, it is known that not only the
standard $n$--sphere, but also exotic spheres of
dimension $\geq 7$ admit such functions. Therefore, 
only by the existence of such a special Morse function one
cannot detect the difference of differentiable structures.

On the other hand, there is a notion of special
generic maps, which constitute a natural generalization of the
above class of Morse functions. Special generic maps
are smooth maps $f : M \to \R^p$ with $1 \leq p \leq n$
that have only definite folds as their singularities
(for details, see Section~\ref{section2}).
As far as the author knows,
this class of maps was first studied by Calabi
\cite{Calabi} who used such a class of maps
in order to study the curvature pinching problem
in differential geometry, although details were not given.
Later, Burlet--de Rham \cite{BdR} studied special generic
maps of $3$--dimensional manifolds into $\R^2$ and
completely characterized those $3$--dimensional closed manifolds
which admit such maps. In \cite{Saeki1}, the
author studied special generic maps on homotopy spheres
and showed that a homotopy $n$--sphere is diffeomorphic
to the standard $n$--sphere if and only if it admits
a special generic map into $\R^p$ for every $1 \leq p \leq n$
(see also \cite{PF}). In \cite{Saeki4}, the author showed
that if we consider cobordisms of Morse functions
with exactly local minima and maxima as their singularities, with
the help of special generic maps into $\R^2$, then we can
recover the $h$--cobordism group of homotopy spheres.
In this way, special generic maps can give us
information on the differentiable structures of manifolds.
It is now known that the existence of special generic maps can, 
in many cases, characterize the standard differentiable
structures (for example, see \cite{Saeki3, SS1, SS2}).

More recently, in \cite{Wr1}, Wrazidlo studied
special generic maps of homotopy spheres with standard
properties (for details, see Section~\ref{section2} of
the present paper), and showed that the
existence of such maps is closely related to
the Gromoll filtration of the $h$--cobordism
group $\Theta_n$ of homotopy $n$--spheres introduced by Gromoll
\cite{Gr1}. More precisely, it was shown that
$\Gamma^n_p \subset F^n_p$ for $n \geq 7$ and $1 \leq p
\leq n-1$, where $\Gamma^n_p$ is the Gromoll
filtration of $\Theta_n$ and $F^n_p$ is the subgroup of 
$\Theta_n$ consisting of homotopy $n$--spheres admitting
standard special generic maps into $\R^p$.

In this paper we sharpen the result by showing that, in fact,
$\Gamma^n_p = F^n_p$ holds. This will be shown by carefully
using the structure theorem for special generic maps
obtained in \cite{Saeki1}.
Note that the Gromoll filtration is still being studied
and the known results are summarized in \cite[Appendix]{CSS}.

The paper is organized as follows. In Section~\ref{section2},
we recall the basic definitions concerning special generic maps
and recall the structure theorem. In Section~\ref{section3},
we recall the definition of the Gromoll filtration \cite{Gr1}.
Several formulations are known in the literature, and
we give one which is useful in the proof of our main result.
Finally, in Section~\ref{section4}, we prove that
$\Gamma^n_p \supset F^n_p$ holds by decomposing the source
homotopy $n$--sphere in a convenient way by using
the structure theorem.
We also give a remark about Calabi's contribution \cite{Calabi}.

Throughout the paper, manifolds and maps
are differentiable of class $C^\infty$
unless otherwise indicated.

\section{Preliminaries}\label{section2}

In this section, we review some known results about
special generic maps which will be necessary for
our purpose. For most of the materials,
the reader is referred to \cite{Saeki1}.

Let $M$ be a closed $n$--dimensional manifold and $p$
an integer with $1 \leq p \leq n$.
For a smooth map $f : M \to \R^p$,
we denote the set of singular points by
$$S(f) = \{x \in M\,|\, \rank{df_x} < p \}.$$

\begin{dfn}\label{dfn1}
(1) A singular point $x \in S(f)$ of $f$ is called a
\emph{definite fold} if there exist
local coordinates $(x_1, x_2, \ldots, x_n)$ 
and $(y_1, y_2, \ldots, y_p)$ around
$x$ and $f(x)$, respectively, such that
$f$ has the form
$y_j \circ f = x_j, j = 1, 2, \ldots, p-1$ and
$y_p \circ f = x_p^2 + x_{p+1}^2 +  \cdots + x_n^2$.

(2) A smooth map $f : M \to \R^p$ is called a \emph{special
generic map} if $S(f)$ consists only of definite folds.
\end{dfn}

It is easy to show that for a special generic map
as above, $S(f)$ is a regular submanifold of $M$
of dimension $p-1$ and that $f|_{S(f)}: S(f)
\to \R^p$ is an immersion.

\begin{exam}\label{ex1}
Consider the unit $n$--sphere $S^n \subset \R^{n+1}$
and the standard projection $pr : \R^{n+1} \to \R^p$, $1 \leq p \leq n$.
Then, one can easily check that the restriction $f = pr|_{S^n} : S^n \to \R^p$
is a special generic map.
\end{exam}

\begin{dfn}
Let $M$ be a closed $n$--dimensional manifold.
The (possibly empty) finite set of integers
$$\mathcal{S}(M) = \{ p \,|\, 1 \leq p \leq n, \, \exists f : M \to
\R^p \mbox{ \rm special generic map}\}$$
is called the \emph{special generic dimensions} of $M$.
This is clearly a diffeomorphism invariant of $M$.
\end{dfn}

By Example~\ref{ex1}, we see that $\mathcal{S}(S^n)$
coincides with the full set $\{1, 2, \ldots, n\}$.

In the following, we mainly consider special generic
maps of $n$--dimensional manifolds into $\R^p$ with $1 \leq p < n$.
For special generic maps into $\R^n$, the reader is referred to \cite{E1}.
For example, it is known that an orientable closed $n$--dimensional
manifold admits a special generic map into $\R^n$ if and only if
$M$ is stably parallelizable. When the manifold $M$ is nonorientable
of dimension $n$, it admits a special generic map into $\R^n$
if and only if there exists a set of $n$ nowhere dependent sections
of $TM \oplus \varepsilon$, where $TM$ is the tangent bundle of $M$
and $\varepsilon$ is the trivial line bundle over $M$
(see \cite{Ando} and \cite[Corollary~2.4]{SSS}).

\begin{dfn}\label{stein}
Let $f: M \to \R^p$ be a special generic map as defined above,
where we assume $1 \leq p < n$.
Two points
$x, x'\in M$ are \emph{$f$--equivalent}
if $f(x)=f(x')$ and the points $x$ and $x'$
are in the same connected component of $f^{-1}(f(x))
= f^{-1}(f(x'))$.
We denote by $W_f$ the quotient space with 
respect to the $f$--equivalence, 
endowed with the quotient topology.
The quotient map is denoted by $q_f: M \to W_f$. 
Then there exists a unique continuous map 
$\bar{f}: W_f \to \R^p$ 
such that $f=\bar{f}\circ q_f$.
The quotient space $W_f$ or the commutative diagram
\begin{eqnarray*}
M \!\!\!\! & \spmapright{f} & \!\!\!\! \R^p \\
& {}_{q_f}\!\!\searrow \quad \qquad \nearrow_{\bar{f}} & \\
& \,W_f &
\end{eqnarray*}
is called the \emph{Stein factorization} of $f$.
The space $W_f$ is often called the \emph{Reeb space}
of $f$ as well.
\end{dfn}

Then, the following structure theorem is known (see \cite{BdR, Saeki1}).

\begin{prop}\label{prop:structure}
Let $f : M \to \R^p$ be a special generic map
of a closed $n$--dimensional manifold $M$ into an Euclidean space
with $1 \leq p < n$.
Then, $W_f$ can be given the unique structure of a smooth $p$--dimensional
compact manifold with nonempty boundary in such a way that $\bar{f} : W_f
\to \R^p$ is an immersion. Furthermore, we have the following.
\begin{enumerate}
\item The restriction $q_f|_{S(f)} : S(f) \to \partial W_f$ is a diffeomorphism.
\item The restriction $q_f|_{M \setminus S(f)} : M \setminus
S(f) \to \Int{W_f}$ is a smooth $S^{n-p}$--bundle.
\item Let $C \cong \partial W_f \times [0, \varepsilon]$, $\varepsilon > 0$,
be a closed collar neighborhood of $\partial W_f$ in $W_f$, $\pi : C \to \partial W_f$
a natural projection, and set $B = q_f^{-1}(C)$. Then, the map $\pi \circ q_f|_B :
B \to \partial W_f$ is a $D^{n-p+1}$--bundle whose structure group is
the orthogonal group $O(n-p+1)$.
\end{enumerate}
\end{prop}

As a corollary, we have the following \cite{Saeki1}.

\begin{cor}\label{cor:structure}
Let $M$ be a closed $n$--dimensional manifold. For $1 \leq p < n$,
there exists a special generic map $f : M \to \R^p$ if and only if
the following holds.
\begin{enumerate}
\item For a compact parallelizable $p$--dimensional
manifold $W$ with nonempty boundary, there exists a smooth
$S^{n-p}$--bundle $E$ over $W$ such that
the structure group of its restriction to $\partial W$ 
can be reduced to the orthogonal group $O(n-p+1)$.
\item There exists a $D^{n-p+1}$--bundle $B$ over $\partial W$ whose
structure group is the orthogonal group $O(n-p+1)$.
\item There exists a bundle isomorphism
$\varphi : \partial B \to \partial E$ between
the smooth $S^{n-p}$--bundles over $\partial W$ such that
$M$ is diffeomorphic to the manifold $E \cup_\varphi B$
obtained by attaching $E$ and $B$ along their boundaries by
the diffeomorphism $\varphi$.
\end{enumerate}
\end{cor}

\begin{rmk}\label{rmk:structure}
In the above corollary, if there exists a special
generic map $f : M \to \R^p$, then the manifold $W$
can be chosen so that it is diffeomorphic to $W_f$.
\end{rmk}

Note that the structure group of the smooth $S^{n-p}$--bundle
$E$ over $W$ as above is the full diffeomorphism group
$\mathrm{Diff}(S^{n-p})$ in general.
When $n-p=1, 2$ or $3$, the natural inclusion
$O(n-p+1) \to \mathrm{Diff}(S^{n-p})$ is known to be
a weak homotopy equivalence (see \cite{Smale1, Hatcher1}).
Furthermore, for the group $\mathrm{Homeo}(D^{n-p+1})$
of homeomorphisms, the restriction map 
$\mathrm{Homeo}(D^{n-p+1}) \to \mathrm{Homeo}(S^{n-p})$
has a right inverse. This implies that every
smooth $S^{n-p}$--bundle over $W$ can be extended to
a topological $D^{n-p+1}$--bundle.
Using these, we have the following \cite{Saeki1,
Saeki2}.

\begin{cor}
Let $M$ be a closed $n$--dimensional manifold.
\begin{enumerate}
\item For $p$ with $n-p = 1, 2$ or $3$, there exists
a special generic map $f : M \to \R^p$ if and only if
$M$ is diffeomorphic to the boundary of a
$D^{n-p+1}$--bundle with structure group $O(n-p+1)$
over a compact parallelizable $p$--dimensional manifold with
nonempty boundary.
\item If there exists a special generic map
$f : M \to \R^p$ for some $1 \leq p < n$,
then $M$ is homeomorphic to the boundary
of a topological $D^{n-p+1}$--bundle over a compact parallelizable
$p$--dimensional manifold with nonempty boundary.
\end{enumerate}
\end{cor}

Now, let us consider special generic maps of homotopy spheres.
First, we recall the following (see \cite{Saeki1}).

\begin{lem}
For a special generic map $f : M \to \R^p$ of a closed
$n$--dimensional manifold $M$ with $1 \leq p < n$,
$M$ is a homotopy $n$--sphere if and only if
the Reeb space $W_f$ is contractible.
\end{lem}

The following definition is originally due to Wrazidlo
\cite{Wr1}.

\begin{dfn}
A special generic map $f : M \to \R^p$ of a closed
$n$--dimensional manifold $M$ with $1 \leq p < n$
is said to be \emph{standard} if the Reeb space $W_f$
is diffeomorphic to the $p$--dimensional disk $D^p$.
\end{dfn}

Note that for $1 \leq p \leq 3$, every special generic map
is standard. For $p=3$, this follows from the positive solution
to the $3$--dimensional Poincar\'e conjecture.
On the other hand, for $p \geq 6$, a special
generic map $f : M \to \R^p$ is standard if and only if
$S(f)$ is simply connected. Furthermore, for all $4 \leq p < n$,
it is known that there exist special generic maps
of homotopy $n$--spheres into $\R^p$ which are not standard.

Let us recall the following lemma \cite[Corollary~3.9]{Wr1}.

\begin{lem}\label{lem:ssgm1}
For a homotopy sphere $\Sigma$, if there exists
a standard special generic map $\Sigma \to \R^p$,
then for each $p'$ with $1 \leq p' < p$, $\Sigma$ has
a standard special generic map into $\R^{p'}$.
\end{lem}

\begin{dfn}
For a homotopy $n$--sphere $\Sigma$, we set
$$\mathcal{F}(\Sigma) = \max\{p\,|\,
1 \leq p < n, \, \exists f : \Sigma \to \R^p \mbox{ \rm
standard special generic map}\},$$
which is called the \emph{fold perfection} of $\Sigma$.
%
\end{dfn}

Note that for $n \geq 5$, we always have
$\mathcal{F}(\Sigma) \geq 2$ by \cite{PF, Saeki1}.
Note also that for an arbitrary homotopy sphere $\Sigma$, we have
$$\mathcal{S}(\Sigma) \supset \{1, 2, \ldots, \mathcal{F}(\Sigma), n\}$$
by Lemma~\ref{lem:ssgm1} and \cite{E1}.

For example, we have
$\mathcal{F}(S^n) = n-1$
by Example~\ref{ex1}.

Now, for $n \geq 5$, let $\Theta_n$ denote
the $h$--cobordism group of oriented homotopy $n$--spheres
\cite{KM}, where the addition corresponds to the
connected sum operation.
For $1 \leq p \leq n-1$, we denote by $F^n_p$
the subset of $\Theta_n$ consisting of (the $h$--cobordism
classes of) homotopy
$n$--spheres $\Sigma$ with $\mathcal{F}(\Sigma) \geq p$.
By virtue of \cite[Lemma~5.4]{Saeki1}, we see that
$F^n_p$ is in fact a subgroup of $\Theta_n$. 
Thus, we have the following filtration of $\Theta_n$ by subgroups:
$$\Theta_n = F^n_1 = F^n_2 \supset \cdots
\supset F^n_{n-4} \supset F^n_{n-3} =
F^n_{n-2} = F^n_{n-1} = 0,$$
where the vanishing of the last three subgroups
follows from \cite[Corollary~4.2]{Saeki1} and
\cite[Remark~2.4]{Saeki2}.

\section{Gromoll filtration}\label{section3}

In this section, we recall the notion of the Gromoll
filtration for the group of homotopy spheres \cite{Gr1}.

For $n \geq 6$, 
it is known that $\Theta_n$ is isomorphic to
the group $\Gamma^n = \pi_0(\mathrm{Diff}^+(S^{n-1}))$ of isotopy classes
of orientation preserving diffeomorphisms
of $S^{n-1}$, where the isomorphism is given by
associating to an orientation preserving diffeomorphism $\varphi : S^{n-1}
\to S^{n-1}$ the homotopy $n$--sphere $D^n \cup_\varphi (-D^n)$
obtained by attaching two copies of the $n$--dimensional disk
along their boundaries (see \cite{C1}).

For $1 \leq p \leq n-1$, let $\mathrm{Diff}^+_p(S^{n-1})$
be the subgroup of $\mathrm{Diff}^+(S^{n-1})$
consisting of diffeomorphisms $\varphi : S^{n-1} \to S^{n-1}$
satisfying $\pi^n_p \circ \varphi = \pi^n_p$, where
$\pi^n_p : S^{n-1} \to \R^{p-1}$
is the standard projection $\R^n \to \R^{p-1}$
restricted to the unit sphere $S^{n-1}$.
The image of the natural map $\pi_0(\mathrm{Diff}^+_p(S^{n-1}))
\to \pi_0(\mathrm{Diff}^+(S^{n-1})) = \Gamma^n$ is denoted
by $\Gamma^n_p$. Thus, we have the following filtration
of $\Gamma^n \cong \Theta_n$ by subgroups:
$$\Gamma^n = \Gamma^n_1 = \Gamma^n_2 \supset
\Gamma^n_3 \supset \cdots \supset \Gamma^n_{n-4}
\supset \Gamma^n_{n-3} = \Gamma^n_{n-2}
= \Gamma^n_{n-1} = 0,$$
where $\Gamma^n_1 = \Gamma^n_2$ follows from Cerf's
pseudo-isotopy theorem \cite{C1}, and
the vanishing of the last three subgroups follows
from \cite{Smale1, Hatcher1}.

Note that the restriction $\pi^n_p : S^{n-1} \to 
\R^{p-1}$ of the natural projection is a special generic map
(see Example~\ref{ex1}).
Therefore, by virtue of Proposition~\ref{prop:structure},
we see that $\pi_0(\mathrm{Diff}^+_p(S^{n-1}))$
is naturally isomorphic to $\pi_{p-1}(\mathrm{Diff}^+(S^{n-p}),
SO(n-p+1))$ (see also \cite{Gr1}).

\section{Fold filtration and Gromoll filtration}\label{section4}

Wrazidlo \cite{Wr1} showed that $\Gamma^n_p \subset
F^n_p$ for $n \geq 7$ and $1 \leq p \leq n-1$.
In fact, we have the following.

\begin{thm}\label{thm1}
For $n \geq 6$ and $1 \leq p \leq n-1$, we have
$\Gamma^n_p = F^n_p$.
\end{thm}

Note that this theorem is trivially valid for
$n=6$, since $\Theta_n = \Gamma^n = 0$.

\begin{proof}[Proof of Theorem~\textup{\ref{thm1}}]
Let us show that $\Gamma^n_p \supset
F^n_p$ for $n \geq 7$ and $1 \leq p \leq n-1$.
Let $\Sigma$ be a homotopy $n$--sphere belonging
to $F^n_p$. Then, there exists a standard special
generic map $f : \Sigma \to \R^p$.
Using the notation as in Corollary~\ref{cor:structure},
we have $\Sigma = E \cup_\varphi B$
and $W$ is diffeomorphic to $D^p$ (see Remark~\ref{rmk:structure}), 
where $E$ is a smooth $S^{n-p}$--bundle over $W$ and $B$
is a $D^{n-p+1}$--bundle over $\partial W \cong S^{p-1}$
with structure group $SO(n-p+1)$, as the latter
is an orientable bundle.
Since $W$ is contractible, $E$ is diffeomorphic
to $S^{n-p} \times W \cong S^{n-p} \times D^p$.

Let us decompose $\partial W \cong S^{p-1}$
into the union of two $(p-1)$--dimensional
disks $\Delta_0$ and $\Delta_1$ attached along their
sphere boundaries.
Over each $\Delta_j$, $j = 0, 1$, $B$
is the trivial $D^{n-p+1}$--bundle, $D^{n-p+1}
\times \Delta_j$.
By modifying the identification diffeomorphism
$E \cong S^{n-p} \times D^p$ if necessary,
we may assume that $\varphi:
S^{n-p} \times \Delta_0 (\subset \partial B) \to 
S^{n-p} \times \Delta_0 (\subset \partial E)$
is the identity map. 

Now, the homotopy $n$--sphere $\Sigma$
is decomposed into the union of 
\begin{eqnarray*}
D_0 & = & E \cup_\varphi (D^{n-p+1} \times \Delta_0)
\cong (S^{n-p} \times W) \cup_{\id}
(D^{n-p+1} \times \Delta_0) \\
& \cong & (S^{n-p} \times D^p) \cup_{\id}
(D^{n-p+1} \times \Delta_0)
\end{eqnarray*}
and $D_1 = D^{n-p+1} \times \Delta_1$ attached along
their boundaries.
Note that $D_0$ and $D_1$ are diffeomorphic to
the $n$--dimensional disk and that $\pi^n_p :
\partial D_1 \cong S^{n-1} \to \R^{p-1}$
corresponds to the restriction of the projection
$D_1 = D^{n-p+1} \times \Delta_1 \to \Delta_1$
to the second factor in an appropriate manner.
In fact, such a correspondence is a special case of 
Corollary~\ref{cor:structure} applied to the canonical
special generic map $\pi^n_p : S^{n-1} \to \R^{p-1}$
as in Example~\ref{ex1}.

Let us now consider the attaching diffeomorphism
$\psi : \partial D_1 \to \partial D_0$,
where $\partial D_1 = (D^{n-p+1} \times \partial \Delta_1)
\cup (S^{n-p} \times \Delta_1)$.
On $D^{n-p+1} \times \partial \Delta_1$,
the diffeomorphism is a bundle map
over $\partial \Delta_1$ with structure group
$SO(n-p+1)$, since the structure group of $B$ 
is given by $SO(n-p+1)$.
Thus, $\psi$ corresponds to an element of
$\pi_{p-1}(\mathrm{Diff}^+(S^{n-p}), SO(n-p+1))
\cong \pi_0(\mathrm{Diff}^+_p(S^{n-1}))$.
This, in turn, corresponds to an element of $\Gamma^n_p$.
Therefore, $\Sigma$ belongs to $\Gamma^n_p$.
This completes the proof.
\end{proof}

For various computations about the Gromoll
filtration $\Gamma^n_p$, the reader is referred to
\cite{CS1} and \cite[Appendix]{CSS}.
For example, the odd multiples of a
generator of $\Gamma^7 \cong
\Z/28\Z$ do not belong to $\Gamma^7_4$ by \cite{We1}.
Therefore, they do not admit a special
generic map into $\R^3$ as pointed out in \cite{Wr1}, although
we still do not know which even multiples of a generator
admit such special generic maps.
 
By \cite{CSS}, we have $\Gamma^n_p \neq 0$
for $(n, p) = (8j+1, 8j-5), (8j+2, 8j-4), (8j+3, 8j-3)$, $j \geq 1$,
and for $(n, p) = (4i-1, 2i-2)$, $i \geq 4$, 
which implies that for these dimensions,
there exist exotic $n$--spheres admitting standard special
generic maps into $\R^p$.

\begin{rmk}
In a survey paper of Calabi 
\cite[Theorem~2]{Calabi},
it is claimed that if $f : \Sigma \to \R^p$
is a special generic map such that $S(f)$
is homeomorphic to $S^{p-1}$ and $f|_{S(f)}$
is an embedding, then $\Sigma$ belongs to
Milnor's group $\Gamma_{n, p-1}$.
However, as far as the author knows, no
proof is given, nor the precise definition
of $\Gamma_{n, p-1}$.
\end{rmk}

We end this paper by posing
a problem.

\begin{prob}
For a homotopy $n$--sphere $\Sigma$,
if there exists a special generic map $f : \Sigma 
\to \R^p$ for some $p$ with $3 < p < n-3$,
then does there exist a standard special
generic map $\Sigma \to \R^p$ ?
\end{prob}

\section*{Acknowledgment}\label{ack}
The author would like to thank
Dominik Wrazidlo for posing the problem,
stimulating discussions and giving important comments
about an earlier version of the paper, which
drastically improved the presentation of the paper.
The author has been supported in part by JSPS KAKENHI Grant Numbers 
JP17H06128, JP22K18267. 
This work was also supported by the Research  
Institute for Mathematical Sciences, an International Joint
Usage/Research Center located in Kyoto University. 
Finally, the author would like to express his hearty thanks
to Professor Maria Aparecida Soares Ruas for her energetic
and constant encouragement.

\end{document}